\documentclass[12pt,reqno]{amsart}
\usepackage{amscd,amsmath,amsthm,amssymb}
\usepackage{color}
\usepackage{pstricks}
\usepackage{stmaryrd}
\usepackage{tikz}
\usepackage{url}

\usepackage{amssymb}
\usepackage{amsfonts,amsmath,mathtools}
\usepackage{graphics}
\usepackage{float}
\usepackage{subfig}
\usepackage{xcolor}
\usepackage{tikz-cd}
\usepackage{booktabs} 
\usepackage{colortbl}
\usepackage[all]{xy}

\usepackage{latexsym}
\usepackage{amsfonts,amsmath,mathtools}
\usepackage{graphics}
\usepackage{float}
\usepackage{enumitem}

\newpsstyle{fatline}{linewidth=1.5pt}
\newpsstyle{fyp}{fillstyle=solid,fillcolor=verylight}
\definecolor{verylight}{gray}{0.97}
\definecolor{light}{gray}{0.9}
\definecolor{medium}{gray}{0.85}
\definecolor{dark}{gray}{0.6}

\def\NZQ{\mathbb}               

\def\ZZ{{\NZQ Z}}

\def\G{{\mathcal G}}

\def\HS{\textup{HS}}
\def\pd{\textup{proj}\phantom{.}\!\textup{dim}}

\def\opn#1#2{\def#1{\operatorname{#2}}} 
\opn\chara{char} \opn\length{\ell} \opn\pd{pd} \opn\rk{rk}
\opn\projdim{proj\,dim} \opn\injdim{inj\,dim} \opn\rank{rank}
\opn\depth{depth} \opn\grade{grade} \opn\height{height}
\opn\embdim{emb\,dim} \opn\codim{codim}

\opn\Tr{Tr} \opn\bigrank{big\,rank}
\opn\superheight{superheight}\opn\lcm{lcm}
\opn\trdeg{tr\,deg}
\opn\reg{reg} \opn\lreg{lreg} \opn\ini{in} \opn\lpd{lpd}
\opn\size{size} \opn\sdepth{sdepth}
\opn\link{link}\opn\fdepth{fdepth}\opn\lex{lex}
\opn\tr{tr}
\opn\type{type}
\opn\gap{gap}
\opn\diam{diam}
\opn\Mod{Mod}
\opn\div{div} \opn\Div{Div} \opn\cl{cl} \opn\Cl{Cl}
\opn\Spec{Spec} \opn\Supp{Supp} \opn\supp{supp} \opn\Sing{Sing}
\opn\Ass{Ass} \opn\Min{Min}\opn\Mon{Mon}
\opn\Ann{Ann} \opn\Rad{Rad} \opn\Soc{Soc}
\opn\Im{Im} \opn\Ker{Ker} \opn\Coker{Coker} \opn\Am{Am}
\opn\Hom{Hom} \opn\Tor{Tor} \opn\Ext{Ext} \opn\End{End}
\opn\Aut{Aut} \opn\id{id}

\opn\nat{nat}
\opn\pff{pf}
\opn\Pf{Pf} \opn\GL{GL} \opn\SL{SL} \opn\mod{mod} \opn\ord{ord}
\opn\Gin{Gin} \opn\Hilb{Hilb}\opn\sort{sort}
\opn\PF{PF}\opn\Ap{Ap}
\opn\dist{dist}
\opn\aff{aff}
\opn\relint{relint} \opn\st{st}
\opn\lk{lk} \opn\cn{cn} \opn\core{core} \opn\vol{vol}  \opn\inp{inp} \opn\nilpot{nilpot}
\opn\link{link} \opn\star{star}\opn\lex{lex}\opn\set{set}
\opn\width{wd}
\opn\Fr{F}
\opn\QF{QF}
\opn\G{G}
\opn\type{type}\opn\res{res}
\opn\conv{conv}
\opn\sr{sr}
\opn\gr{gr}

\def\pot#1#2{#1[\kern-0.28ex[#2]\kern-0.28ex]}

\opn\dirlim{\underrightarrow{\lim}}
\opn\inivlim{\underleftarrow{\lim}}
\def\Implies{\ifmmode\Longrightarrow \else
	\unskip${}\Longrightarrow{}$\ignorespaces\fi}
\def\implies{\ifmmode\Rightarrow \else
	\unskip${}\Rightarrow{}$\ignorespaces\fi}
\def\iff{\ifmmode\Longleftrightarrow \else
	\unskip${}\Longleftrightarrow{}$\ignorespaces\fi}

\let\:=\colon
\newtheorem{Theorem}{Theorem}[section]
\newtheorem{Lemma}[Theorem]{Lemma}
\newtheorem{Corollary}[Theorem]{Corollary}

\newtheorem{Example}[Theorem]{Example}

\newtheorem{Definition}[Theorem]{Definition}

\newtheorem{Conjecture}[Theorem]{Conjecture}
\newtheorem{Characterization}[Theorem]{Characterization}

\newtheorem{Criterion}[Theorem]{Criterion}

\let\epsilon\varepsilon
\let\kappa=\varkappa
\def\qed{\ifhmode\textqed\fi
	\ifmmode\ifinner\hfill\quad\qedsymbol\else\dispqed\fi\fi}
\def\textqed{\unskip\nobreak\penalty50
	\hskip2em\hbox{}\nobreak\hfill\qedsymbol
	\parfillskip=0pt \finalhyphendemerits=0}
\def\dispqed{\rlap{\qquad\qedsymbol}}

\opn\dis{dis}
\def\pnt{{\raise0.5mm\hbox{\large\bf.}}}

\opn\Lex{Lex}

\let\emptyset\varnothing
\let\epsilon\varepsilon

\def\Ass{\textup{Ass}}

\def\HS{\textup{HS}}
\def\lcm{\textup{lcm}}
\def\gcd{\textup{gcd}}
\def\set{\textup{set}}

\begin{document}
	
\title{Very well--covered graphs via the Rees algebra}
\author{Marilena Crupi, Antonino Ficarra}

\address{Marilena Crupi, Department of mathematics and computer sciences, physics and earth sciences, University of Messina, Viale Ferdinando Stagno d'Alcontres 31, 98166 Messina, Italy}
\email{mcrupi@unime.it}

\address{Antonino Ficarra, Department of mathematics and computer sciences, physics and earth sciences, University of Messina, Viale Ferdinando Stagno d'Alcontres 31, 98166 Messina, Italy}
\email{antficarra@unime.it}



\maketitle

\begin{abstract}
A very well--covered graph is a well--covered graph without isolated vertices such that the size of its minimal vertex covers is half of the number of vertices. If $G$ is a Cohen--Macaulay very well--covered graph, we deeply investigate some algebraic properties of the cover ideal of $G$ via the Rees algebra associated to the ideal, and especially when $G$ is a whisker graph.
\end{abstract}

\section{Introduction}

In this article, a \textit{graph} will always mean a finite undirected graph without loops or multiple edges. Let $G$ be a graph with the \textit{vertex set} $V(G)=\{x_1,\dots,x_n\}$ and the \textit{edge set} $E(G)$. Let $K$ be a field and let $R=K[x_1,\dots,x_n]$ be the standard graded polynomial ring with coefficients in $K$. If $x_i\in V(G)$, the \textit{neighborhood} of $x_i$ is the set
\[N(x_i)=\{x_j\in V(G):x_ix_j\in E(G)\}.\]

By abuse of notation, we use an edge $e=\{x_i, x_j\}$ interchangeably with the monomial $x_ix_j\in R$.\smallskip
 
Let $W \subseteq V(G)$. $W$ is called a \emph{vertex cover} if $e\cap W\ne\emptyset$ for all $e\in E(G)$, and it is called a \emph{minimal vertex cover} if no proper subset of $W$ is a vertex cover of $G$. The set of all minimal vertex covers of $G$ is denoted by $\mathcal{C}(G)$. Attached to $G$ \cite{RV1, RV} there are the \emph{edge ideal} of $G$, defined as
\[I(G) = (x_ix_j\ :\ x_ix_j \in E(G)),\]
and the \emph{cover ideal} of $G$, defined as
\[J(G) = (x_{i_1}\cdots x_{i_s}\ :\ \{x_{i_1},\ldots, x_{i_s}\}\in\mathcal{C}(G)).
\]

The cover ideal $J(G)$ is the \emph{Alexander dual} of the edge ideal $I(G)$ and conversely. 
Indeed, the minimal primes of the edge ideals, 
correspond to the minimal vertex covers of their underlying graph.

Hereafter, for an integer $n\ge1$, we set $[n]=\{1,\dots,n\}$.

A graph $G$ is \emph{well--covered} or \emph{unmixed} if all its minimal vertex covers have the same cardinality. In particular, all the associated primes of $I(G)$ have the same height \cite{RV}. By \cite[Corollary 3.4]{GV}, for a well--covered graph $G$ without isolated vertices we have
\[2\cdot\min\{|C|:C\in\mathcal{C}(G)\}\ \ge\ \vert V(G)\vert.\] 
If equality holds, the graph $G$ is called \textit{very well--covered}.
Note that, if $G$ is a very well--covered graph, then its number of vertices is even.

Finally, a graph $G$ is called \emph{Cohen–Macaulay} over the field $K$ if $R/I(G)$ is a Cohen--Macaulay ring, $R=K[x_i: x_i\in V(G)]$. It is well--known that a Cohen--Macaulay graph is well--covered \cite[Proposition 7.2.9]{RV}. The fundamental Eagon--Reiner theorem \cite[Theorem 8.1.9]{JT} says that $R/I(G)$ is a Cohen--Macaulay ring if and only if $J(G)$ has a linear resolution.

Very well--covered graphs have been studied from view points of both Commutative Algebra and Combinatorics. See, for instance, \cite{CRT2011,OF,GV,KTY2018,KPFTY,KPTY2022,MMCRTY2011,Fak}. 
In \cite{CF2023} we deeply studied very well--covered graphs by means of \textit{Betti splittings} \cite{FHT2009}. Recently, many authors have managed the Betti splitting technique for studying algebraic and combinatorial properties of classes of monomial ideals (see, for instance, \cite{CFts1,CFL,FH2023} and references therein).

In the present article, we continue the algebraic study of the class of Cohen--Macaulay very well--covered graphs started in \cite{CF2023}. If $G$ is a graph in such a class, our main tool will be the Rees algebra of the cover ideal $J(G)$. We state that if $G$ is a Cohen--Macaulay very well--covered graph, then the Rees algebra of $J(G)$ is a normal Cohen--Macaulay domain and as a consequence we obtain some relevant properties on the behavior of the powers of $J(G)$, when $G$ is a \emph{whisker graph}. Adding a \emph{whisker} to a graph $G$ at a vertex $v$ means adding a new vertex $w$ and an edge $vw$ to the set $E(G)$. If a whisker is added to every vertex of $G$, then the resulting graph, denoted by $G^\ast$, is called the \emph{whisker graph} or \emph{suspension} of $G$. It is important to point out that the whisker graph $G^\ast$ of a graph $G$ with $n$ vertices is a very well-covered Cohen--Macaulay graph with $2n$ vertices (see, for instance, \cite{KPFTY} and references therein). In \cite[Question 6.6]{HVT2023} it is asked in which way attaching whiskers to a graph $H$ gives rise to a graph $G$ such that $J(G)$ has linear powers. In Corollary \ref{Conj:PowersHSCMverywell} we partially answer this question. See also \cite[Corollary 4.5]{LW2023} and \cite[Theorem 2.3]{MF2014}.

Here are the outline of the article. In Section \ref{sec:2} we discuss a normality criterion for squarefree monomial ideals (Criterion \ref{Thm:CritNormality}). This result is borrowed from \cite{NQBM}. Section \ref{sec:3} deeply investigates the Rees algebra $\mathcal{R}(J(G))$ of $J(G)$, with $G$ a Cohen--Macaulay very well--covered graph. Our main result states that $\mathcal{R}(J(G))$ is a normal Cohen--Macaulay domain (Theorem \ref{Thm:R(J(G))CMnormal}). To obtain this result we use Criterion \ref{Thm:CritNormality} as well as the structure theorem of Cohen--Macaulay very well--covered graphs (Characterization \ref{char:veryWellCGCM}) stated in \cite{CRT2011}. In Section \ref{sec:4}, if $G$ is a whisker graph with $2n$ vertices, we prove that $J(G)$ satisfies the $\ell$-exchange property (Theorem \ref{Thm:J(G)ell-exchangeProp}). As a consequence, we state that $J(G)^k$ has linear quotients, for all $k\ge 1$, and then that $J(G)$ has linear powers (Corollary \ref{Cor:J(G)LinQuot}). Furthermore, if $G$ is a whisker graph, we show that each power of $J(G)$ has homological linear quotients (Theorem \ref{Thm:HSJ(G)LinQuotRees}). This result supports a conjecture stated in \cite[Conjecture 4.4]{CF2023}. Moreover, as applications of the previous results, we compute the limit depth, the depth stability and the analytic spread of $J(G)$. Finally,  if  $G$ is a whisker graph with $2n$ vertices, we get a partial result on the structure of the reduced Gr\"obner basis of the presentation ideal of $\mathcal{R}(J(G))$ (Corollary \ref{cor:grobner}).  At present we do not know the reduced Gr\"obner basis of the presentation ideal of $\mathcal{R}(J(G))$. However, our experiments in \textit{Macaulay2} \cite{GDS} suggest that for a suitable monomial order, the reduced Gr\"obner basis is quadratic and hence that $\mathcal{R}(J(G))$ is Koszul (Conjecture \ref{Conj:R(J(G))Koszul}), for any Cohen--Macaulay very well--covered graph with $2n$ vertices.

\section{A normality criterion for monomial ideals}\label{sec:2}

Let $I$ be an ideal of a domain $R$. An element $f\in R$ is \textit{integral over $I$} if it satisfies an equation of the type
$$
f^k+a_1f^{k-1}+\dots+a_{k-1}f+a_k=0, \ \ \ a_i\in I^i.
$$
The set of all these elements, denoted by $\overline{I}$, is an ideal containing $I$ and called the \textit{integral closure} of $I$. We say that $I$ is \textit{integrally closed} if $\overline{I}=I$, and we say that $I$ is \textit{normal} if all its powers $I^k$, $k\ge1$, are integrally closed.

Let $I$ be an ideal of a commutative ring $R$ generated by $u_1, \ldots, u_m$. The \textit{Rees algebra} of $I$, denoted by $\mathcal{R}(I)$ or $R[It]$, is the subring of $R[t]$, defined as follows
$$
\mathcal{R}(I)= R[It]= R[u_1t, \ldots, u_mt]= \bigoplus_{k\ge0}I^kt^k \subset\ R[t],
$$
where $t$ is a new variable. We quote the next fundamental result from \cite{HV} (see, also, \cite[Theorem 4.3.17]{RV}).

\begin{Theorem} \label{thm:HVnormality} Let $I$ be an ideal of a normal domain $R$. Then the
following are equivalent:
\begin{enumerate}
\item[\em(a)] $I$ is a normal ideal;
\item[\em(b)] the Rees algebra $\mathcal{R}(I)$ is normal.
\end{enumerate}
\end{Theorem}

This property will be crucial in the sequel.\smallskip

Now, let $R=K[x_1,\dots,x_n]$ be the standard graded polynomial ring with coefficients in a field $K$ and let $I$ be a monomial ideal of $R$. 
As usual we denote by $G(I)=\{u_1,\dots,u_m\}$ the unique minimal set of monomial generators of $I$. Then the Rees algebra of $I$ is the following $K$-algebra
$$
\mathcal{R}(I)\ =\ K[x_1,\dots,x_n,u_1t,\dots,u_mt]\ \subset\ R[t].
$$

The next criterion quickly follows from \cite[Theorem 3.1]{NQBM} (see, also, \cite[Theorem 3.1]{NBR2022}).
\begin{Criterion}\label{Thm:CritNormality}
	Let $I_1,I_2\subset K[x_2,\dots,x_{n}]$ be two squarefree monomial ideals. If $I_1\subseteq I_2$ are normal ideals, then $I=I_1+x_1I_2\subset R$ is a normal squarefree monomial ideal.
\end{Criterion}

	Indeed, by \cite[Theorem 3.1]{NQBM}, it is enough to check that $I_1+I_2$ is normal and that $\gcd(x_1,u)=1$ for all $u\in G(I_1)\cup G(I_2)$. Since $I_1\subseteq I_2$, the first assertion follows because $I_1+I_2=I_2$ is normal by hypothesis. The second  assertion follows because the generators of $I_1$ and $I_2$ are monomials of $K[x_2,\dots,x_n]$.

\section{The Rees algebra}\label{sec:3}
In this section we study the Rees algebra of the vertex cover ideal of a Cohen--Macaulay very well--covered graph. 

%

The main result in this section is the following.
\begin{Theorem}\label{Thm:R(J(G))CMnormal}
	Let $G$ be a Cohen--Macaulay very well--covered graph. Then the Rees algebra $\mathcal{R}(J(G))$ is a normal Cohen--Macaulay domain.
\end{Theorem}



In order to prove it, 
we recall the following fundamental algebraic characterization of Cohen--Macaulay very well--covered graphs. 

\begin{Characterization}\label{char:veryWellCGCM}
	\textup{(\cite{CRT2011}, \cite[Lemma 3.1]{MMCRTY2011}).} Let $G$ be a very well--covered graph with $2n$ vertices. Then, the following conditions are equivalent.
	\begin{enumerate}[label=\textup{(\alph*)}]
		\item $G$ is Cohen--Macaulay.
		\item There exists a labeling of $V(G)=\{x_1,\dots,x_n,y_1,\dots,y_n\}$ such that
		\begin{enumerate}[label=\textup{(\roman*)}]
			\item $X=\{x_1,\dots,x_n\}$ is a minimal vertex cover of $G$ and $Y=\{y_1,\dots,y_n\}$ is a maximal independent set of $G$,
			\item $x_iy_i\in E(G)$ for all $i\in[n]$,
			\item if $x_iy_j\in E(G)$ then $i\le j$,
			\item if $x_iy_j\in E(G)$ then $x_ix_j\notin E(G)$,
			\item if $z_ix_j,y_jx_k\in E(G)$ then $z_ix_k\in E(G)$ for any distinct $i,j,k$ and $z_i\in\{x_i,y_i\}$.
		\end{enumerate}
	\end{enumerate}

\end{Characterization}

From now on, if $G$ is a Cohen--Macaulay very well--covered graph with $2n$ vertices, we tacitly assume that its set of vertices $V(G)=\{x_1,\dots,x_n,y_1,\dots,y_n\}$ satisfies  the conditions \textup{(i)-(v)} of Characterization \ref{char:veryWellCGCM}, without having to relabel it. See \cite{CRT2011, CF2023} for more details on this topic. 

Hereafter, denote by $S$ the polynomial ring $K[x_1,\dots,x_n,y_1,\dots,y_n]$ in the $2n$ variables $x_1,\dots,x_n,y_1,\dots,y_n$ with coefficients in the field $K$.\smallskip

Let $F\subseteq[n]$ be a non empty set. We set ${\bf x}_F=\prod_{i\in F}x_i$, ${\bf y}_F=\prod_{i\in F}y_i$. Otherwise, we set ${\bf x}_{\emptyset}={\bf y}_{\emptyset}=1$. The \textit{support} of a monomial $u\in S$ is the set 
$$
\supp(u)=\{x_i:x_i\ \textit{divides}\ u\}\cup\{y_j:y_j\ \textit{divides}\ u\}.
$$

If $W \subseteq V(G)$, we denote by $G\setminus W$ the subgraph of $G$ with the vertices of $W$ and their incident edges deleted. 

The following results were proved in \cite[Lemma 2.2, Proposition 2.1]{CF2023}.
\begin{Lemma}\label{Lem:minimalGeneratorsI(G)^vee}\label{Prop:GremoveVerticesVeryWC}
	Let $G$ be a Cohen--Macaulay very well--covered graph with $2n$ vertices.
	\begin{enumerate}
		\item[\textup{(a)}] For each $u\in G(J(G))$ there exists a unique subset $F$ of $[n]$ such that $u={\bf x}_F{\bf y}_{[n]\setminus F}$.
		\item[\textup{(b)}] For any $A\subseteq[n]$, $G\setminus\{x_i,y_i:i\in A\}$ is a Cohen--Macaulay very well--covered graph.
	\end{enumerate}
\end{Lemma}

For a subset $C$ of $X\cup Y=\{x_1,\dots,x_n,y_1,\dots,y_n\}$, we define
$$
{\bf z}_C= {\bf x}_{C_x}{\bf y}_{C_y},
$$
where $C_x=\{i:x_i\in C\}$ and $C_y=\{j:y_j\in C\}$. The next result holds true. 
\begin{Lemma}\label{Lemma:J(G)Decompx_1}
	Let $G$ be a Cohen--Macaulay very well--covered graph with $2n$ vertices. Then
	\begin{equation}\label{eq:J(G)x_1Decomp}
	J(G)\ =\ {\bf z}_{N(x_1)}J(G_1)+x_1J(G\setminus\{x_1,y_1\}),
	\end{equation}
	where $G_1=G\setminus\{x_{i},y_{i}:i\in N(x_1)_x\cup N(x_1)_y\}$.
\end{Lemma}
\begin{proof}
	The proof is similar to that of \cite[Proposition 2.3]{CF2023}. We include it for completeness. Let $u\in G(J(G))$. By Lemma \ref{Lem:minimalGeneratorsI(G)^vee}(a), either $x_1$ divides $u$ or $y_1$ divides $u$.\medskip\\
	\textsc{Case 1.} Suppose $x_1$ divides $u$. Note that $N(y_1)=\{x_1\}$. Indeed, by Characterization \ref{char:veryWellCGCM}(i), $N(y_1)$ is a subset of $X$, since $Y$ is a maximal independent set. Moreover, by (iii) if $x_iy_1\in E(G)$ then $i\le 1$. Hence, $N(y_1)=\{x_1\}$. Consequently ${\bf z}_{N(y_1)}=x_1$ and the support $C'$ of $u/{\bf z}_{N(y_1)}=u/x_1$ is a vertex cover of $G\setminus\{x_1,y_1\}$. But $C'$ is a minimal vertex cover, for $u/x_1$ has degree $n-1$ and $G\setminus\{x_1,y_1\}$ is a Cohen--Macaulay very well--covered graph with $2(n-1)$ vertices (Lemma \ref{Prop:GremoveVerticesVeryWC}(b)). Thus $u/x_1\in G(J(G\setminus\{x_1,y_1\}))$ and so $u\in G(x_1J(G\setminus\{x_1,y_1\}))$.\medskip\\
	\textsc{Case 2.} Suppose $y_1$ divides $u$. Since the support $C$ of $u$ is a minimal vertex cover of $G$ and $x_1\notin C$, then $z_i\in C$ for all $z_i\in N(x_1)$. Consequently, the support $C_1$ of $u/{\bf z}_{N(x_1)}$ is a vertex cover of $G_1$. But $C_1$ is a minimal vertex cover of $G_1$, for $|C_1|=n-|N(x_1)|$ and $G_1$ is a Cohen--Macaulay very well--covered graph with $2(n-|N(x_1)|)$ vertices (Lemma \ref{Prop:GremoveVerticesVeryWC}(b)). Hence $u\in G({\bf z}_{N(x_1)}J(G_1))$.\medskip
	
	These two cases show the inclusion ``$\subseteq$" in equation (\ref{eq:J(G)x_1Decomp}). The other inclusion is acquired as in the last part of the proof of \cite[Proposition 2.3]{CF2023}.
\end{proof}

\begin{Corollary}\label{Cor:normal}
	Let $G$ be a Cohen--Macaulay very well--covered graph with $2n$ vertices. Then $J(G)$ is a normal ideal.
\end{Corollary}
\begin{proof}
	By Lemma \ref{Lemma:J(G)Decompx_1}, equation (\ref{eq:J(G)x_1Decomp}) holds. Set $J=J(G)$, $J_1={\bf z}_{N(x_1)\setminus y_1}J(G_1)$ and $J_2=J(G\setminus\{x_1,y_1\})$. Thus
	$$
	J\ =\ y_1J_1+x_1J_2.
	$$
	Since $y_1J_1,J_2\subset K[x_2,\dots,x_n,y_1,y_2,\dots,y_n]$, it is enough to show that $J_1\subseteq J_2$. Then $y_1J_1\subset J_2$ and the result follows from Criterion \ref{Thm:CritNormality} and induction on $n$.
	
	Let $u\in G(J_1)$. We must prove that $u\in G(J_2)$, too. That is, we must show that $C=\supp(u)$ is a minimal vertex cover of $G\setminus\{x_1,y_1\}$. It is enough to prove $C$ is a vertex cover of $G\setminus\{x_1,y_1\}$. Minimality follows because $|C|=n-1$. Hence, we must prove that $e\cap C\ne\emptyset$ for all edges $e\in E(G\setminus\{x_1,y_1\})$. Let $e\in E(G\setminus\{x_1,y_1\})$. Since $y_1u\in G(J)$, it follows that $C\cup y_1$ is a minimal vertex cover of $G$. Hence $e\cap (C\cup y_1)\ne\emptyset$. Therefore $e\cap C\ne\emptyset$ because $y_1\notin e$. Our assertion follows.
\end{proof}

Finally, we are in the position to prove  the main result in the section.
\begin{proof}[Proof of Theorem \ref{Thm:R(J(G))CMnormal}]
By Corollary \ref{Cor:normal}, $J(G)$ is a normal ideal. Hence, the Rees algebra $\mathcal{R}(J(G))$ is normal (Theorem \ref{thm:HVnormality}). Next, by a theorem of Hochster \cite{Hoc72}, since $\mathcal{R}(J(G))$ is a normal affine semigroup ring, it follows that $\mathcal{R}(J(G))$ is Cohen--Macaulay.
\end{proof}

The \textit{toric ring} of $G(J(G))$ is the $K$-algebra $K[J(G)]=K[u:u\in G(J(G))]\subset S$.

\begin{Corollary}
	Let $G$ be a Cohen--Macaulay very well--covered graph. Then the toric ring $K[J(G)]$ is a normal Cohen--Macaulay domain.
\end{Corollary}
\begin{proof}
	Since $J(G)$ is generated in one degree, the statement follows from Theorem \ref{Thm:R(J(G))CMnormal} together with \cite[Proposition 4.3.42]{RV}.
\end{proof}

Now let $I$ be an ideal of a noetherian ring $R$. As usual, denote by $V(I)$ the set of prime ideals containing $I$ and by $\Ass(I)$ the set of associated prime ideals of $R/I$. For all $P\in \Spec(R)$, we denote by $\mathfrak{m}_P$ the maximal ideal of the local ring $R_P$. Recall that $I$ satisfies the \textit{persistence property} (with respect to  associated ideals) if
$$
\Ass(I)\subseteq\Ass(I^2)\subseteq\Ass(I^3)\subseteq\cdots.
$$
In \cite{HQ15}, Herzog and Qureshi introduced the notion of \textit{strong persistence property}. More in detail, let $P\in V(I)$. We say that $I$ satisfies the \textit{strong persistence property with respect to $P$} if for all $k$ and all $f\in (I^k_P:\mathfrak{m}_P) \setminus I^k_P$ there exists $g\in I_P$ such that $fg\notin I^{k+1}_P$. The ideal $I$ is said to satisfy the \textit{strong persistence property} if it satisfies the strong persistence property for all $P\in V(I)$. One can verify that the strong persistence property implies the persistence property \cite{HQ15} (see, also, \cite[Proposition 2.1]{NQBM}).

Theorem \ref{Thm:R(J(G))CMnormal} yields the next result.

\begin{Corollary}
	Let $G$ be a Cohen--Macaulay very well--covered graph. Then $J(G)$ satisfies the strong persistence property, and in particular, the persistence property.
\end{Corollary}
\begin{proof}
	The assertion follows from Theorem \ref{Thm:R(J(G))CMnormal} and \cite[Corollary 1.6]{HQ15}.
\end{proof}

\section{Whisker graphs}\label{sec:4}
In this section we study some algebraic properties of the powers of the cover ideals of a special class of Cohen--Macaulay very well--covered graphs. Our main tool is the so called $\ell$-exchange property introduced in \cite{HHV2005}. 

Let $I\subset R=K[x_1,\dots,x_n]$ be a monomial ideal generated in one degree, and let $K[I]=K[u:u\in G(I)]$ be the toric ring of $G(I)$. Then $K[I]$ has the presentation $\psi:T=K[t_u:u\in G(I)]\rightarrow K[I]$ defined by $\psi(t_u)=u$ for all $u\in G(I)$. The kernel $\Ker(\psi)=J$ is called the \textit{toric ideal} of $K[I]$. 

Fix a monomial order $>$ on $T$. We say that the monomial $t_{u_1}\cdots t_{u_N}\in T$ is \textit{standard with respect to $>$}, if $t_{u_1}\cdots t_{u_N}$ does not belong to the initial ideal, $\text{in}_<(J)$, of the toric ideal $J$ of $K[I]$.

Let $u\in S$ be a monomial. The \textit{$x_i$-degree} and the \textit{$y_i$-degree} of $u$ are the integers $\deg_{x_i}(u)=\max\{j:x_i^j\ \textit{divides}\ u\}$, $\deg_{y_i}(u)=\max\{j:y_i^j\ \textit{divides}\ u\}$, respectively.

\begin{Definition}{\rm(}\cite[Definition 3.3]{DHQ}{\rm)}. 
	The equigenerated monomial ideal $I\subset R$ satisfies the \textit{$\ell$-exchange property with respect to $>$}, if the following condition is satisfied: for all standard monomials $t_{u_1}\cdots t_{u_N}$, $t_{v_1}\cdots t_{v_N}\in T$ of degree $N$ such that
	\begin{enumerate}
		\item[\em(i)] $\deg_{x_i}(u_1\cdots u_N)=\deg_{x_i}(v_1\cdots v_N)$, for all $1\le i\le j-1$ with $j\le n-1$,
		\item[\em(ii)] $\deg_{x_j}(u_1\cdots u_N)<\deg_{x_j}(v_1\cdots v_N)$,
	\end{enumerate}
	there exist $h$ and $k$ with $j< h \le n$ and $1\le k\le N$, such that $x_j(u_k/x_h)\in G(I)$.
\end{Definition}
The following lemmata will be needed later.

\begin{Lemma}\label{Lem:againCover}
	Let $G$ be a Cohen--Macaulay very well--covered graph with $2n$ vertices. Let $C\in\mathcal{C}(G)$ such that $C_y\ne\emptyset$ and let $i=\min C_y$. Then $(C\setminus y_i)\cup x_i\in\mathcal{C}(G)$.
\end{Lemma}
\begin{proof}
	Firstly we prove that $C'=(C\setminus y_i)\cup x_i$ is a vertex cover of $G$. Let $e\in E(G)$, we must show that $e\cap C'$ is non empty. Since $C$ is a vertex cover of $G$, then $e\cap C\ne\emptyset$. If $\{x_i,y_i\}\cap e=\emptyset$, then $e\cap C'\ne\emptyset$, too. If $x_i\in e$ then $e\cap C'$ contains $x_i$ and therefore the intersection is non empty. Finally, suppose $y_i\in e$ but $x_i\notin e$. Since $Y$ is a maximal independent set, it follows that $N(y_i)\subseteq X$. Hence, $e=x_jy_i$ for some $j$. By Characterization \ref{char:veryWellCGCM}(iii) we have $j\le i$. Thus $j<i$, because $x_i\notin e$. Since $i=\min C_y$ and $j<i$, it follows from Lemma \ref{Lem:minimalGeneratorsI(G)^vee}(a) that $x_j\in C$. Hence $x_j\in e\cap C'$ and again the intersection is non empty.
	
	The fact that $C'$ is a minimal vertex cover of $G$ follows because $|C'|=n$.
\end{proof}

\begin{Lemma}\label{Lemma:degx+degy=const}
	Let $G$ be a Cohen--Macaulay very well--covered graph with $2n$ vertices. Then, for all $k\ge1$, all $i\in[n]$ and all $u\in G(J(G)^k)$ we have
	$$
	\deg_{x_i}(u)+\deg_{y_i}(u)=k.
	$$
\end{Lemma}
\begin{proof}
	By Lemma \ref{Lem:minimalGeneratorsI(G)^vee}(a), for all $u\in G(J(G))$, we have
	$$
	\deg_{x_i}(u)+\deg_{y_i}(u)=1\ \ \text{for all}\ \ 1\le i\le n.
	$$
	Since $J(G)$ is generated in a single degree, the minimal generators of $J(G)^k$ are the products $u=u_1\cdots u_k$ of $k$ arbitrary monomials of $G(J(G))$. Hence, for all $1\le i\le n$, we have $\deg_{x_i}(u)+\deg_{y_i}(u)=\sum_{j=1}^k[\deg_{x_i}(u_j)+\deg_{y_i}(u_j)]=k$.
\end{proof}

Now, we consider a special but wide class of Cohen--Macaulay very well--covered graphs. 

Let $H$ be a graph on the vertex set $X=\{x_1,\dots,x_n\}$ and take a new set of variables $Y=\{y_1,\dots,y_n\}$. Then, the \textit{whisker graph} $G=H^*$ of $H$ is the graph obtained from $H$
by attaching to each vertex $x_i$ a new vertex $y_i$ and the edge $x_iy_i$. The edge $x_iy_i$ is called a \textit{whisker}. More in detail, the \textit{whisker graph} $G=H^*$ of $H$ is the graph on the vertex set $X\cup Y=\{x_1,\dots,x_n\}\cup\{y_1,\dots,y_n\}$ and the edge set $E(G)\cup\{x_1y_1,x_2y_2,\dots,x_ny_n\}$.

We have already underlined in the introduction that the whisker graph of a given graph with $n$ vertices is a Cohen--Macaulay very well--covered graph with $2n$ vertices (see, also, \cite[Corollary 4.3]{CF2023}, proof). 

With the same notation as before, the next result holds true.

\begin{Lemma}\label{Lem:againCoverWhisker}
	Let $G$ be a whisker graph with vertex set $X\cup Y$. 
	Then, 
	for all $C\in\mathcal{C}(G)$ and all $y_i\in C$ we have that $(C\setminus y_i)\cup x_i\in\mathcal{C}(G)$.
\end{Lemma}
\begin{proof}
	By assumption $G=H^*$, for some graph $H$ on the vertex set $X= \{x_1, \ldots, x_n\}$, thus a Cohen--Macaulay very well--covered graph with the $2n$ vertices $x_1,\dots,x_n,$ $y_1,\dots,$ $y_n$.
Since  for all $i$, the only vertex adjacent to $y_i$ is $x_i$, then for any labeling of $X\cup Y$, the conditions (i)--(v) of Characterization \ref{char:veryWellCGCM} are satisfied. Hence, if $C\in\mathcal{C}(G)$ and $y_i\in C$, we can choose a labeling such that $\min_yC=i$. The assertion follows by applying Lemma \ref{Lem:againCover}.
\end{proof}

From now on, when we tell about a whisker graph $G$ with $2n$ vertices, we implicitly assume that 
\begin{enumerate}
\item[-] $G$ is the whisker graph associated to a given graph whose vertex set is the set $X=\{x_1,\dots,x_n\}$, and with whiskers $x_iy_i$, $i=1, \ldots, n$, that is, $V(G)=X\cup Y$, with $Y=\{y_1,\dots,y_n\}$. 
\item[-] $G$ is a very well--covered Cohen--Macaulay graph whose vertex set $X\cup Y$ satisfies the conditions  (i)--(v) of Characterization \ref{char:veryWellCGCM}. 
\end{enumerate}


\begin{Theorem}\label{Thm:J(G)ell-exchangeProp}
	 Let $G$ be a whisker graph with $2n$ vertices. 
	 Then $J(G)$ satisfies the $\ell$-exchange property with respect  to the lexicographic order $>_{\lex}$ induced by $x_1>y_1>x_2>y_2>\cdots>x_n>y_n$.
\end{Theorem}
\begin{proof}
	Set $z_{2p-1}=x_{p}$ and $z_{2p}=y_{p}$, for $p=1,\dots,n$. Then
	$$
	z_{1}>z_{2}>z_{3}>z_{4}>\dots>z_{2n-1}>z_{2n}.
	$$
	
	We prove the following slightly more general statement.\smallskip\\
	$(*)$ For all monomials $t_{u_1}\cdots t_{u_N}$, $t_{v_1}\cdots t_{v_N}$ of $K[t_u:u\in G(J(G))]$ such that
	\begin{enumerate}
		\item[(i)] $\deg_{z_i}(u_1\cdots u_N)=\deg_{z_i}(v_1\cdots v_N)$, for all $1\le i\le j-1$ with $j\le 2n-1$,
		\item[(ii)] $\deg_{z_j}(u_1\cdots u_N)<\deg_{z_j}(v_1\cdots v_N)$,
	\end{enumerate}
	there exist $h$ and $k$, with $j<h\le 2n$ and $1\le k\le N$, such that $z_j(u_k/z_h)\in G(J(G))$.\medskip
	
	Let $t_{u_1}\cdots t_{u_N}$, $t_{v_1}\cdots t_{v_N}$ monomials of $K[t_u:u\in G(J(G))]$ satisfying the conditions (i) and (ii).
	
	We claim that the integer $j$ is odd. Suppose for a contradiction that $j$ is even, then $j=2p$ for some $p\in[n]$. Thus $z_j=y_p$ and
	\begin{equation}\label{eq:degxy1}
		\deg_{y_p}(u_1\cdots u_N)<\deg_{y_p}(v_1\cdots v_N).
	\end{equation}
	On the other hand, since $J(G)$ is generated in a single degree, $u_1\cdots u_N$ and $v_1\cdots v_N$ belong to $G(J(G)^N)$. Thus, Lemma \ref{Lemma:degx+degy=const} gives
	\begin{align}
		\label{eq:degxy2}\deg_{x_p}(u_1\cdots u_N)+\deg_{y_p}(u_1\cdots u_N)=\deg_{x_p}(v_1\cdots v_N)+\deg_{y_p}(v_1\cdots v_N)=N.
	\end{align} 
	Equations (\ref{eq:degxy1}) and (\ref{eq:degxy2}) yield $\deg_{x_p}(u_1\cdots u_N)>\deg_{x_p}(v_1\cdots v_N)$, but this contradicts condition (i), since $x_p=z_{j-1}$. Hence $j$ is odd, and so $z_j=x_{p}$ for some $p\in[n]$.
	
	Since $\deg_{x_p}(u_1\cdots u_N)<\deg_{x_p}(v_1\cdots v_N)\le N$, by Lemma \ref{Lemma:degx+degy=const} it follows that $\deg_{y_p}(u_1\cdots u_N)>0$. Hence, there exists $k$ with $1\le k\le N$ such that $y_p$ divides $u_k$. By Lemma \ref{Lem:againCoverWhisker}, it follows that $x_p(u_k/y_p)\in G(J(G))$. Since $y_p=z_{j+1}$, and $z_j>z_{j+1}$, the claim $(*)$ is proved.
\end{proof}

Recall that an ideal $I$ of a polynomial ring $R=K[x_1,\dots,x_n]$ is said to have \textit{linear powers} if $I^k$ has linear resolution, for all $k\ge1$. Moreover, a monomial ideal $I$ of $R$, has \textit{linear quotients} if for some order $u_1,\dots,u_m$ of its minimal generating set $G(I)$, all colon ideals $(u_1,\dots,u_{\ell-1}):u_{\ell}$, $\ell=2,\dots,m$, are generated by a subset of the set of variables $\{x_1,\dots,x_n\}$.

As a first consequence of Theorem \ref{Thm:J(G)ell-exchangeProp} we prove 
	that 
	the cover ideal of a whisker graph $G$ has linear powers. Such a result has been recently obtained in \cite[Corollary 4.5]{LW2023} 
(see, also, \cite[Theorem 2.3]{MF2014}) 
by showing that the ordinary powers of 
$J(G)$ are weakly polymatroidal, which implies having linear powers.

\begin{Corollary}\label{Cor:J(G)LinQuot}
	Let $G$ be a whisker graph with $2n$ vertices. Then,
	\begin{enumerate}
	\item[\em(a)] for all $k\ge1$, $J(G)^k$ has linear quotients with respect  to the lexicographic order $>_{\lex}$ induced by $x_1>y_1>x_2>y_2>\cdots>x_n>y_n$.
	\item[\em(b)] $J(G)$ has linear powers. In particular, the depth function $\depth S/J(G)^k$ is a non-increasing function of $k$, that is, $\depth S/J(G)^k\ge\depth S/J(G)^{k+1}$ for all $k\ge1$.
	\end{enumerate}
\end{Corollary}
\begin{proof}
	(a) Since $J(G)$ is generated in a single degree, each minimal monomial generator of $J(G)^N$ is a product $u_1\cdots u_N$ of $N$ arbitrary, non necessarily distinct, monomials $u_i\in G(J(G))$. Let $u=u_1\cdots u_N\in G(J(G)^N)$, where each $u_i\in G(J(G))$. Setting $P=(v_1\cdots v_N:v_i\in G(J(G)),v_1\cdots v_N>_{\lex}u_1\cdots u_N)$, we must prove that the ideal $P:u$ is generated by variables.
	
	Let $v=v_1\cdots v_N\in G(P)$. Using the labeling $z_i$ on the variables, given in the proof of Theorem \ref{Thm:J(G)ell-exchangeProp}, by the definition of $>_{\lex}$, for some $i$ and $j$ we have
	\begin{enumerate}
		\item[(i)] $\deg_{z_i}(v_1\cdots v_N)=\deg_{z_i}(u_1\cdots u_N)$, for all $1\le i\le j-1$ with $j\le 2n-1$,
		\item[(ii)] $\deg_{z_j}(v_1\cdots v_N)>\deg_{z_j}(u_1\cdots u_N)$.
	\end{enumerate}
	Hence, by the property $(*)$ proved in Theorem \ref{Thm:J(G)ell-exchangeProp}, there exist integers $k$ and $h$ such that $z_j=x_{h}$ and $x_{h}(u_k/y_h)\in G(J(G))$. Since $x_h>y_h$, we have $x_{h}(u_k/y_{h})>_{\lex}u_k$. Consequently,
	$$
	u'=x_h(u/y_h)=u_1\cdots u_{k-1}\cdot x_h(u_k/y_h)\cdot u_{k+1}\cdots u_N\in P
	$$
	and $x_h=u'/\gcd(u',u)\in P:u$ divides the monomial $v/\gcd(v,u)\in P$.  Indeed, the set $\{v/\gcd(v,u):v\in G(P)\}$ generates $P:u$ (\cite[Proposition 1.2.2]{JT}). Hence, we see that $P:u$ is generated by variables, as desired.\\
	(b) That $J(G)$ has linear powers follows from (a) 
	and the fact that all powers $J(G)^k$ are monomial ideals generated in a single degree. The claim about the non-increasingness of the function $\depth S/J(G)^k$ follows from \cite[Proposition 10.3.4]{JT}.
\end{proof}

Let us briefly recall the concept of \textit{homological shift ideal} \cite{Bay019,BJT019,CF2023,F2,F2Pack,FH2023,HMRZ021a,HMRZ021b}. For ${\bf a}=(a_1,\dots,a_n)\in\ZZ_{\ge0}^n$, we set ${\bf x^a}=x_1^{a_1}\cdots x_n^{a_n}$. Let $I\subset S$ be a monomial ideal, then $\HS_i(I)=({\bf x^a}:\beta_{i,\bf a}(I)\ne0)$ is the \emph{$i$th homological shift ideal} of $I$, where $\beta_{i,\bf a}(I)$ is a multigraded Betti number. Note that $\HS_0(I)=I$ and $\HS_i(I)=0$ if $i<0$ or $i>\pd(I)$. A basic goal is to determine those homological properties satisfied by all $\HS_i(I)$, the so-called \textit{homological shift properties} of $I$ \cite{F2}. In \cite{CF2023}, we posed the following conjecture:
\begin{Conjecture}\label{Conj:PowersHSCMverywell}
	\textup{(\cite[Conjecture 4.4]{CF2023}).} Let $G$ be a Cohen--Macaulay very well--covered graph with 2n vertices. Then $\HS_k(J(G)^\ell)$ has linear quotients, for all $k\ge0$, and all $\ell\ge1$.
\end{Conjecture}
In \cite{CF2023}, we gave a positive answer to this conjecture for $\ell=1$ \cite[Theorem 4.4]{CF2023} and for all Cohen--Macaulay bipartite graphs \cite[Corollary 4.11]{CF2023}. Now we prove that the powers of cover ideals of whisker graphs have homological linear quotients, partially answering Conjecture \ref{Conj:PowersHSCMverywell}.
	\begin{Theorem}\label{Thm:HSJ(G)LinQuotRees}
		Let $G$ be a whisker graph with $2n$ vertices. 
		Then, for all $\ell\ge1$ and all $k\ge0$, $\HS_k(J(G)^\ell)$ has linear quotients with respect to the lexicographic order $>_{\lex}$ induced by $x_1>x_2>\cdots>x_n>y_1>y_2>\cdots>y_n$.
	\end{Theorem}
\begin{proof}
	Let $>$ be the lexicographic order induced by $x_1>y_1>x_2>y_2>\cdots>x_n>y_n$. Then, by Corollary \ref{Cor:J(G)LinQuot}(a), $J(G)^\ell$ has linear quotients with respect to $>$ for all $\ell\ge1$. Let $u\in G(J(G)^\ell)$, we define
	$$
	\set(u)\ =\ \{i\ :\ z_i\in\{x_i,y_i\},z_i\in(v\in G(J(G)^\ell):v>u):(u)\}.
	$$
	The definition of the order $>$ and Lemma \ref{Lemma:degx+degy=const} imply that the set of variables generating the ideal $(v\in G(J(G)^\ell):v>u):(u)$ is a subset of $X=\{x_1,\dots,x_n\}$.
	
	Thus, by \cite[Proposition 1.2]{F2} we have
	$$
	\HS_k(J(G)^\ell)\ =\ \big({\bf x}_Fu\ :\ u\in G(J(G)^{\ell}),\ F\subseteq\set(u),\ |F|=k\big).
	$$
	
	Let ${\bf x}_Dv\in G(\HS_k(J(G)^{\ell}))$, $D\subseteq\set(v)$, $v\in G(J(G)^{\ell})$ and consider the colon ideal
	$$
	P\ =\ ({\bf x}_Fu\in G(\HS_k(J(G)^{\ell})):{\bf x}_Fu>_{\lex}{\bf x}_Dv):({\bf x}_Dv).
	$$
	We must prove that $P$ is generated by variables.

	Let ${\bf x}_Fu\in G(\HS_k(J(G)^{\ell}))$, $F\subseteq\set(u)$, $u\in G(J(G)^\ell)$, such that ${\bf x}_Fu>_{\lex}{\bf x}_Dv$. Let $h=\lcm({\bf x}_Fu,{\bf x}_Dv)/({\bf x}_Dv)$. If $\deg(h)=1$, $h$ is a variable. Assume $\deg(h)>1$. Let $z_i$ be the labeling on the variables such that $z_1=x_1$, $z_2=x_2$, $\dots$, $z_n=x_n$, $z_{n+1}=y_1$, $z_{n+2}=y_2$, $\dots$, $z_{2n}=y_n$. Then, by definition of $>_{\lex}$, there exists $p$ such that $\deg_{z_j}({\bf x}_Fu)=\deg_{z_j}({\bf x}_Dv)$ for all $j<p$ and
	\begin{equation}\label{eq:ineqdeg}
		\deg_{z_p}({\bf x}_Fu)>\deg_{z_p}({\bf x}_Dv).
	\end{equation}
	
	Now, we distinguish two cases.\smallskip\\
	\textsc{Case 1.} Suppose $z_p=x_i$ for some $i$. We claim that
	\begin{equation}\label{eq:boundell}
		\deg_{x_i}({\bf x}_Fu)\le\ell.
	\end{equation}
	Indeed, by Lemma \ref{Lemma:degx+degy=const} and the structure of ${\bf x}_Fu$, it follows that $\deg_{x_i}({\bf x}_Fu)\le\ell+1$. Suppose by contradiction that $\deg_{x_i}({\bf x}_Fu)=\ell+1$, then $i\in\set(u)$. Necessarily $y_i$ must divide $u$. But this would imply that $\deg_{x_i}({\bf x}_Fu)+\deg_{y_i}({\bf x}_Fu)$ exceeds $\ell+1$, which is impossible. Hence, equation (\ref{eq:boundell}) follows.
	
	By Lemma \ref{Lemma:degx+degy=const} and equations (\ref{eq:ineqdeg}) and (\ref{eq:boundell}), we have $\ell\ge\deg_{x_i}({\bf x}_Fu)>\deg_{x_i}({\bf x}_Dv)$ and $\deg_{y_i}({\bf x}_Dv)>0$. Writing $v=v_1v_2\cdots v_\ell$, with each $v_q\in G(J(G))$, we have that $y_i$ divides $v_q$ for some $q$. Then $i\in\set(v)$. Indeed, $x_i(v_q/y_i)\in G(J(G))$ and $v'=x_i(v/y_i)=v_1\cdots v_{q-1}(x_i(v_q/y_i))v_{q+1}\cdots v_{\ell}\in G(J(G)^\ell)$. We distinguish two cases.\smallskip\\
	\textsc{Subcase 1.1.} Let $i\notin D$. Then ${\bf x}_Dv'\in G(\HS_k(J(G)^\ell))$ and ${\bf x}_Dv'>_{\lex}{\bf x}_Dv$. Moreover, $\lcm({\bf x}_Dv',{\bf x}_Dv)/({\bf x}_Dv)=x_i\in P$ divides $h$.\smallskip\\
	\textsc{Subcase 1.2.} Let $i\in D$. Then $\deg_{x_i}({\bf x}_Dv)+\deg_{y_i}({\bf x}_Dv)=\ell+1$ (Lemma \ref{Lemma:degx+degy=const}). Since by (\ref{eq:ineqdeg}) and (\ref{eq:boundell}) $\deg_{x_i}({\bf x}_Dv)<\ell$, it follows that $\deg_{y_i}(v)\ge2$. Therefore, there exist $h_1\ne h_2$ such that $y_i$ divides $v_{h_1}$ and $v_{h_2}$. Set $v'=v_1\cdots (x_i(v_{h_1})/y_i)\cdots v_{h_2}\cdots v_\ell$. Then, it follows that $v'\in G(J(G)^\ell)$ and $D\subseteq\set(v')$. Thus ${\bf x}_Dv'>_{\lex}{\bf x}_Dv$ and moreover $\lcm({\bf x}_Dv',{\bf x}_Dv)/({\bf x}_Dv)=x_i\in P$ divides $h$.\medskip\\
	\textsc{Case 2.} Suppose $z_p=y_i$ for some $i$. For all $j$ such that $\deg_{y_j}({\bf x}_Fu)>\deg_{y_j}({\bf x}_Dv)$, since $\deg_{x_j}({\bf x}_Fu)=\deg_{x_j}({\bf x}_Dv)$, Lemma \ref{Lemma:degx+degy=const} gives
	\begin{align*}
		\ell+1=\deg_{x_j}({\bf x}_Fu)+\deg_{y_j}({\bf x}_Fu)>\deg_{x_j}({\bf x}_Dv)+\deg_{y_j}({\bf x}_Dv)=\ell.
	\end{align*}
    Hence $j\notin D$ and $\deg_{y_j}({\bf x}_Fu)-\deg_{y_j}({\bf x}_Dv)=1$. Let $j_1,\dots,j_t$ be the integers such that $\deg_{y_{j_s}}({\bf x}_Fu)>\deg_{y_{j_s}}({\bf x}_Dv)$, $s=1,\dots,t$. Then, the above argument shows that $h=y_{j_1}y_{j_2}\cdots y_{j_t}$ and $j_s\notin D$ for all $s=1,\dots,t$. Since $\deg(h)>1$ we have $t\ge2$. As before $v'=x_{j_2}\cdots x_{j_t}v/(y_{j_2}\cdots y_{j_t})\in G(J(G)^{\ell})$, $D\subseteq\set(v')$ and ${\bf x}_Dv'>_{\lex}{\bf x}_Dv$. Finally $\lcm({\bf x}_Dv',{\bf x}_Dv)/({\bf x}_Dv)=y_{j_1}\in P$ divides $h$ .\medskip
    
    The above \textsc{Cases 1} and 2 show that $P$ is generated by variables, as wanted.
\end{proof}

Another relevant consequence of Theorem \ref{Thm:J(G)ell-exchangeProp} concerns the \textit{limit depth} of $J(G)$. The role of the Rees algebra of the cover ideal $J(G)$ will be crucial  to calculating it.\medskip

Let $I$ be a graded ideal of a polynomial ring $R$ with $n$ variables. By a theorem of Brodmann \cite{Brod79}, $\depth R/I^k$ is constant for $k$ large enough. This eventually constant value is called the \textit{limit depth} of $I$, and it is denoted by $\lim_{k\rightarrow\infty}\depth R/I^k$.

The \textit{depth stability} of $I$, denoted by $\textup{dstab}(I)$, is the least integer $k_0$ such that $\depth R/I^k=\depth R/I^{k_0}$ for all $k\ge k_0$. Brodmann proved that
$$
\lim_{k\rightarrow\infty}\depth R/I^k\le n-\ell(I),
$$
where $\ell(I)$ is the \textit{analytic spread} of $I$, that is, the Krull dimension of the fiber ring $\mathcal{R}(I)/\mathfrak{m}\mathcal{R}(I)$, where $\mathfrak{m}$ is the maximal graded ideal of $R$. 

If the Rees algebra of $I$ is Cohen--Macaulay, then by \cite[Proposition 10.3.2]{JT} (see, also, \cite{EH83} combined with \cite[Proposition 1.1]{Huneke82}), we have
\begin{equation}\label{eq:limdepthCM}
\lim_{k\rightarrow\infty}\depth R/I^k=n-\ell(I).
\end{equation}

Hence, we have
\begin{Theorem}\label{Thm:ell(J(G))R(J(G))}
	Let $G$ be a whisker graph with $2n$ vertices. 
	Then
	$$
	\lim_{k\rightarrow\infty}\depth S/J(G)^k\ =\ n-1.
	$$
	Moreover, $\textup{dstab}(J(G))\le n$ and $\ell(J(G))=n+1$.
\end{Theorem}

For the proof of this result, we need the next more general lemma. Let $u\in S$ be a monomial. Using the notation in Section \ref{sec:3}, we have $\supp(u)_x=\{i:x_i\ \textit{divides}\ u\}$ and $\supp(u)_y=\{i:y_i\ \textit{divides}\ u\}$.


\begin{Lemma}\label{Lem:min(y)=i}
	Let $G$ be a Cohen--Macaulay very well--covered graph with $2n$ vertices. Then, for all $i\in[n]$, there exists $u\in G(J(G))$ such that $\min\supp(u)_y=i$.
\end{Lemma}
\begin{proof}
	We proceed by induction on $n\ge1$. For the base case, $J(G)=(x_1,y_1)$ and the statement holds. Now, let $n>1$. Then, by Lemma \ref{Lemma:J(G)Decompx_1} we have $J(G)={\bf z}_{N(x_1)}J(G_1)+x_1J(G\setminus\{x_1,y_1\})$. If $u\in G({\bf z}_{N(x_1)}J(G_1))$, then $\min\supp(u)_y=1$. Let $i\in[n]$ with $i\ne1$. Since $G\setminus\{x_1,y_1\}$ is a Cohen--Macaulay very well--covered graph (Lemma \ref{Lem:minimalGeneratorsI(G)^vee}(b)), by induction there exists $v\in G(J(G\setminus\{x_1,y_1\}))$ with $\min\supp(v)_y=i$. Setting $u=x_1v$, we have $u\in G(J(G))$ and $\min\supp(u)_y=\min\supp(v)_y=i$. The proof is complete.
\end{proof}
\begin{proof}[Proof of Theorem \ref{Thm:ell(J(G))R(J(G))}]
	For any $N\ge1$ and any $u_1\cdots u_N\in G(J(G)^N)$, with $u_i\in G(J(G))$, let $P=(v_1\cdots v_N:v_i\in G(J(G)),v_1\cdots v_N>_{\lex}u_1\cdots u_N)$ and denote by $s(u_1,\dots,u_N)$ the number of variables generating the colon ideal $P:u_1\cdots u_N.$
	Since $J(G)^N$ has linear quotients with respect to $>_{\lex}$ (Corollary \ref{Cor:J(G)LinQuot}(a)) we have that $\depth S/J(G)^{N}=2n-\max\{s(u_1,\dots,u_N)+1:u_1,\dots,u_N\in G(J(G))\}$. This follows from \cite[Corollary 8.2.2]{JT} and the Auslander--Buchsbaum formula. The definition of the order $>_{\lex}$ and Lemma \ref{Lemma:degx+degy=const} imply that the set of variables generating the ideal $P:u_1\cdots u_N$ is a subset of $X=\{x_1,\dots,x_n\}$. Pick $n$ monomials $u_i={\bf x}_{F_i}{\bf y}_{[n]\setminus F_i}$, where $\min([n]\setminus F_i)=i$, for $i=1,\dots,n$. The existence of these monomials follows from Lemma \ref{Lem:min(y)=i}. By Lemma \ref{Lem:againCover}, $x_i(u_i/y_i)\in G(J(G))$ for $i=1,\dots,n$. Hence, $s(u_1,\dots,u_n)=n$. This shows that $\depth S/J(G)^n=n-1$. We claim that $\depth S/J(G)^N=n-1$ for all $N\ge n$. It is enough to consider $u_1,\dots,u_n$ and $N-n$ arbitrary monomials $v_{n+1},\dots,v_{N}\in G(J(G))$. Then $s(u_1,\dots,u_n,v_{n+1},\dots,v_n)=n$ and $\depth S/J(G)^N=n-1$. Hence, $\textup{dstab}(J(G))\le n$.  Moreover, from Theorem \ref{Thm:R(J(G))CMnormal} and equation (\ref{eq:limdepthCM}), since $S$ is a polynomial ring in $2n$ variables, $\ell(J(G))=2n -\lim_{k\rightarrow\infty}\depth S/J(G)^k = n+1$.
\end{proof}

We close the section with some remarks on the reduced Gr\"obner basis of the presentation ideal $J$ of $\mathcal{R}(J(G))$.

Hereafter, we follow closely \cite[Section 6.4.1]{EHGB}.\smallskip

Let $G$ be a Cohen--Macaulay very well--covered graph with $2n$ vertices $X\cup Y=\{x_1,\dots,x_n,y_1,\dots,y_n\}$. Let $G(J(G))=\{u_1,\dots,u_m\}$ and let $\mathcal{R}(J(G))$ be the Rees algebra of $J(G)$. Let ${\bf x}=x_1,\dots,x_n$, ${\bf y}=y_1,\dots,y_n$ and ${\bf t}=t_{u_1},\dots,t_{u_m}$.

Then the Rees algebra $\mathcal{R}(J(G))$ has the presentation
$$
\varphi: S' =K[{\bf x},{\bf y},{\bf t}]\rightarrow\mathcal{R}(J(G))
$$
defined by setting $\varphi(x_i)=x_i$, $\varphi(y_i)=y_i$, for $1\le i\le n$, and $\varphi(t_{u_j})=u_jt$ for $1\le j\le m$. The ideal $J=\Ker(\varphi)$ is called the \textit{presentation ideal} of $\mathcal{R}(J(G))$.\smallskip

Analogously, the toric ring $K[J(G)]=K[u_1,\dots,u_m]$ has the presentation
$$
\psi:T=K[{\bf t}]\rightarrow K[u_1,\dots,u_m]
$$
defined by setting $\psi(t_{u_j})=u_j$ for $1\le j\le m$. The ideal $L=\Ker(\psi)$ is called the \textit{toric ideal} of $K[u_1,\dots,u_m]$. 
Let $\mathfrak{m}=(x_1,\dots,x_n,y_1,\dots,y_n)$ be the graded maximal ideal of $S$. Since $J(G)$ is generated in a single degree, the fiber ring $\mathcal{R}(J(G))/\mathfrak{m}\mathcal{R}(J(G))$ is isomorphic to the toric ring $K[J(G)]$.\\

Let $>'$ be an arbitrary monomial order on $T$ and let $>_{\lex}$ be the lexicographic order on $S$ induced by $x_1>y_1>x_2>y_2>\dots>x_n>y_n$. We define the monomial order $>'_{\lex}$ as follows: for two monomials $w_1t_{u_1}^{a_1}\cdots t_{u_m}^{a_m}$ and $w_2t_{u_1}^{b_1}\cdots t_{u_m}^{b_m}$ in $S'$, with $w_1,w_2\in S$, we set $w_1t_{u_1}^{a_1}\cdots t_{u_m}^{a_m}>_{\lex}'w_2t_{u_1}^{b_1}\cdots t_{u_m}^{b_m}$ if and only if $w_1>_{\lex}w_2$ or $w_1=w_2$ and $t_{u_1}^{a_1}\cdots t_{u_m}^{a_m}>'t_{u_1}^{b_1}\cdots t_{u_m}^{b_m}$. According to \cite[Section 2]{EHGB} the order $>'_{\lex}$ is the product order of $>'$ and $>_{\lex}$.\bigskip

With the notation above, from Theorem \ref{Thm:J(G)ell-exchangeProp} and \cite[Theorem 6.24]{EHGB}, we get the next result. 

\begin{Corollary}\label{cor:grobner}
	Let $G$ be a whisker graph with $2n$ vertices. 
	Then the reduced Gr\"obner basis of the presentation ideal $J$ of $\mathcal{R}(J(G))$ with respect to $>_{\lex}'$ consists of all binomials belonging to the reduced Gr\"obner basis of $L$ with respect to $>'$ together with the binomials
	$$
	x_it_{u}-y_it_{x_i(u/y_i)},
	$$
	where $u,x_i(u/y_i)\in G(J(G))$.
\end{Corollary}

The statement of Corollary \ref{cor:grobner} seems to be true for all Cohen--Macaulay very well--covered graphs.
\begin{Example}
	\rm 
	Consider the graph $G$ with 12 vertices depicted below
	
	
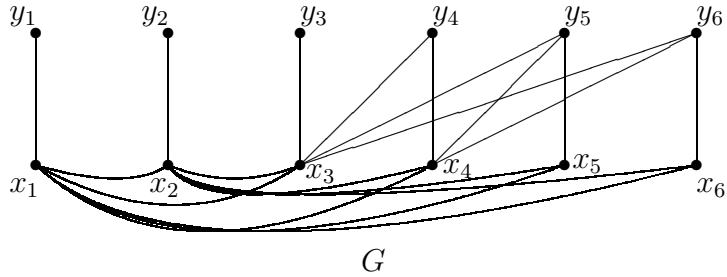
\begin{picture}(90,150)(-90,-50)
		\put(-10,60){\circle*{4}}
		\put(40,60){\circle*{4}}
		\put(-20,65){\textit{$y_1$}}
		\put(30,65){\textit{$y_2$}}
		\put(-10,10){\circle*{4}}
		\put(40,10){\circle*{4}}
		\put(-20,0){\textit{$x_1$}}
		\put(33,0){\textit{$x_2$}}
		\put(90,10){\circle*{4}}
		\put(92,5){\textit{$x_3$}}
		\put(90,60){\circle*{4}}
		\put(90,65){\textit{$y_3$}} 
		\put(-10,60){\line(0,-1){50}}
		\put(40,60){\line(0,-1){50}}
		\put(90,60){\line(0,-1){50}} 
		\put(90,10){\line(1,1){50}}
		\put(144,7){\textit{$x_4$}}
		\put(140,10){\circle*{4}}
		\put(140,65){\textit{$y_4$}}
		\put(140,60){\circle*{4}}
		\put(193,9){\textit{$x_5$}} 
		\put(190,10){\circle*{4}}
		\put(190,65){\textit{$y_5$}}
		\put(190,60){\circle*{4}} 
		\put(190,60){\line(0,-1){50}}
		\put(240,0){\textit{$x_6$}} 
		\put(240,10){\circle*{4}}
		\put(240,65){\textit{$y_6$}}
		\put(240,60){\circle*{4}} 
		\put(240,60){\line(0,-1){50}}
		\put(140,60){\line(0,-1){50}}
		\put(190,60){\line(-2,-1){100}} 
		\put(190,60){\line(-1,-1){50}} 
		\put(240,60){\line(-2,-1){100}} 
		\put(240,60){\line(-3,-1){150}} 
		\qbezier(-10,10)(25,0)(40,10)
		\qbezier(-10,10)(45,-20)(90,10)
		\qbezier(-10,10)(45,-40)(140,10) 
		\qbezier(-10,10)(45,-40)(190,10) 
		\qbezier(-10,10)(45,-40)(240,10) 
		\qbezier(40,10)(65,0)(90,10)
		\qbezier(40,10)(55,-12)(140,10) 
		\qbezier(40,10)(55,-12)(190,10) 
		\qbezier(40,10)(55,-12)(240,10) 
		\put(113,-30){\textit{$G$}}
	\end{picture}\vspace*{-0,3cm}
	
	By Characterization \ref{char:veryWellCGCM}, one verifies that $G$ is a Cohen--Macaulay very well--covered graph with $12$ vertices. We have 
	\begin{align*}
		I(G)\ &\ =\  (x_1y_1, x_2y_2, x_3y_3, x_4y_4, x_5y_5, x_6y_6, x_1x_2, x_1x_3, x_1x_4, x_1x_5,\\
		&\phantom{\ =(..}x_1x_6, x_2x_3, x_2x_4, x_2x_5, x_2x_6, x_3y_4,x_3y_5,x_3y_6, x_4y_5, x_4y_6),\\[4pt]
		J(G)\ &\ =\ (x_1x_2x_3x_4x_5x_6,\, x_1x_2x_3x_4x_5y_6,\, x_1x_2x_3x_4y_5x_6,\, x_1x_2x_3x_4y_5y_6,   \\
		&\phantom{\ =(..} x_1x_2x_3y_4y_5y_6,\, x_1x_2y_3y_4y_5y_6,\, x_1y_2x_3x_4x_5x_6,\, y_1x_2x_3x_4x_5x_6).
	\end{align*}

	We order the monomials $u_1,\dots,u_8$ of $G(J(G))$ with respect to the lexicographic order induced by $x_1>y_1>\dots>x_6>y_6$. Thus, for instance $u_1=x_1x_2x_3x_4x_5x_6$, $u_2=x_1x_2x_3x_4x_5y_6$ and so on.
	
	Now, let
	$$
	\varphi: S' =K[{\bf x},{\bf y},{\bf t}]\rightarrow\mathcal{R}(J(G))
	$$
	be the map defined by setting $\varphi(x_i)=x_i$, $\varphi(y_i)=y_i$, for $1\le i\le 6$, and $\varphi(t_{u_j})=u_jt$ for $1\le j\le 8$. Furthermore, let $T=K[{\bf t}]$.
	
	Let $>'_{\lex}$ be the product order of the lexicographic order $>_{\lex}$ on $S$ induced by $x_1>y_1>\dots>x_n>y_n$, and the lexicographic order $>'$ on $T$ such that $t_{u_i}>'t_{u_j}$ if and only if $u_i>_{\lex}u_j$. By using \textit{Macaulay2} \cite{GDS}, we have that the reduced Gr\"obner basis of $\ker(\varphi)$ with respect to the order $>'_{\lex}$ is the following one:
	\begin{align*}
	\mathcal{G}\ =\ \{&
		x_{6}t_{u_2}-y_{6}t_{u_1},\,\,x_{5}t_{u_3}-y_{5}t_{u_1},\,\,x_{6}t_{u_4}-y_{6}t_{u_3},\,\,x_{5}t_{u_3}-y_{5}t_{u_1},\,\,x_{4}t_{u_5}-y_{4}t_{u_4},\,\,
		\\&x_{3}t_{u_6}-y_{3}t_{u_5},\,\,x_{2}t_{u_7}-y_{2}t_{u_1},\,\,x_{1}t_{u_8}-y_{1}t_{u_1},\,\,t_{u_1}t_{u_4}-t_{u_2}t_{u_3}\}.
	\end{align*}
\end{Example}\medskip

Our experiments using \textit{Macaulay2} \cite{GDS} suggest the next conjecture.

\begin{Conjecture}\label{Conj:R(J(G))Koszul}
Let $G$ be a Cohen--Macaulay very well--covered graph with $2n$ vertices. Then the presentation ideal of the Rees algebra of $J(G)$ has a quadratic reduced Gr\"obner basis with respect to the product order $>_{\lex}'$ of the lexicographic order $>_{\lex}$ on $S$ induced by $x_1>y_1>\dots>x_n>y_n$, and the lexicographic order $>'$ on $T$ such that $t_{u_i}>'t_{u_j}$ if and only if $u_i>_{\lex}u_j$. In particular, $\mathcal{R}(J(G))$ is a Koszul algebra.
\end{Conjecture}

\emph{Acknowledgment}. We thank S.A. Seyed Fakhari and M. Nasernejad for their comments and helpful suggestions that allowed us to improve the quality of the paper.

\end{document}